\documentclass[12pt]{article}

\usepackage[latin1]{inputenc}

\usepackage{enumerate}
\usepackage{amsthm}
\usepackage{a4}
\usepackage{amsmath}

\usepackage{makeidx} 

\usepackage{graphicx} 
\usepackage{color} 
\usepackage{multirow}
\usepackage{amssymb}

\usepackage[colorlinks=true,urlcolor=blue,citecolor=blue,linkcolor=blue,bookmarks=true,bookmarksopen=true,pdfstartview=FitH]{hyperref}
\newcommand {\junk}[1]{}

\theoremstyle{plain}
\newtheorem{th-def}{Theorem-Definition}[section]
\newtheorem{theo}[th-def]{Theorem}
\newtheorem{lem}[th-def]{Lemma}
\newtheorem{prop}[th-def]{Proposition}
\newtheorem{coro}[th-def]{Corollary}

\theoremstyle{definition}
\newtheorem{defi}[th-def]{Définition}

\theoremstyle{remark}
\newtheorem{rem}{\noindent Remark}[section]

\newtheorem{ex}{\noindent Example}

\def\N{\mathbb N}
\def\Z{\mathbb Z}

\def\R{\mathbb R}

\def\P{\mathbb P}
\def\L{\mathbb L}

\def\E{\mathbb E}
\def\A{\mathcal A}

\def\rma{{\mathbb R}_{\max}}

\def\G{\mathcal{G}}

\def\sAn{\left(A(n)\right)_{n\in\N}}

\def\1{\mathbf{1}}

\title{Law of Large Numbers for products of random matrices with coefficients in the max-plus semi-ring.\\
Version of \today}

\author{\href{mailto:glenn.merlet@gmail.com}{Glenn MERLET}}

\begin{document}
\maketitle

\begin{abstract}
We analyze the asymptotic behavior of  random variables $x(n,x_0)$ defined by $x(0,x_0)=x_0$ and $x(n+1,x_0)=A(n)x(n,x_0)$, where $\sAn$ is a stationary and ergodic sequence of random matrices with entries in the semi-ring \mbox{$\R\cup\{-\infty\}$} whose addition is the $\max$ and whose multiplication is $+$.

Such sequences modelize a large class of discrete event systems, among which timed event graphs, 1-bounded Petri nets, some queuing networks, train or computer networks.
We give necessary conditions for $\left(\frac{1}{n}x(n,x_0)\right)_{n\in\N}$ to converge almost surely. Then, we prove a general scheme to give partial converse theorems. When $\max_{A_{ij}(0)\neq -\infty}|A_{ij}(0)|$ is integrable, it allows us:
\begin{enumerate}[-]
\item to give a necessary and sufficient condition for the convergence of $\left(\frac{1}{n}x(n,0)\right)_{n\in\N}$ when the sequence $\sAn$ is i.i.d.,
\item to prove that, if $\sAn$ satisfy a condition of reinforced ergodicity and a condition of fixed structure (i.e. $\P\left(A_{ij}(0)=-\infty\right)\in\{0,1\}$), then  $\left(\frac{1}{n}x(n,0)\right)_{n\in\N}$ converges almost-surely,
\item and to reprove the convergence of $\left(\frac{1}{n}x(n,0)\right)_{n\in\N}$ if the diagonal entries are never $-\infty$.
\end{enumerate}
\end{abstract}

\section{Introduction}
\subsection{Model}
We analyze the asymptotic behavior of random variables $x(n,x_0)$ defined by:
\begin{equation}\label{eqdefx}
\left\{\begin{array}{lcl}
x(0,x_0)&=&x_0\\
x(n+1,x_0)&=&A(n)x(n,x_0),
\end{array} \right.
\end{equation}
where $\sAn$ is a stationary and ergodic sequence of random matrices with entries in the semi-ring \mbox{$\R\cup\{-\infty\}$} whose addition is the $\max$ and whose multiplication is $+$.

We also define the product of matrices $A^n:=A(n-1) A(n-2)\cdots A(0)$  such that $x(n,x_0)=A^n x_0$ and, if the sequence has indices in $\Z$, which is possible up to a change of probability space, $A^{-n}:=A(-1)\cdots A(-n)$ and $y(n,x_0):=A^{-n}x_0$.

On the coefficients, Relation~(\ref{eqdefx}) reads
$$x_i(n+1,x_0)=\max_j\left(A_{ij}(n)+x_j(n,x_0)\right),$$
and the product of matrices is defined by
\begin{equation}\label{eqdefpdt}
(A^n)_{ij}=\max_{i_0=j,i_n=i}\sum_{l=0}^{n-1}A_{i_{l+1}i_l}(l).
\end{equation}
In most cases, we assume that $A(n)$ never has a line of $-\infty$, which is a necessary and sufficient condition for $x(n,x_0)$ to be finite. (Otherwise, some coefficients can be $-\infty$.)\\

Such sequences modelize a large class of discrete event systems. This class includes some models of operations research like timed event graphs (F.~Baccelli~\cite{Baccelli}), 1-bounded Petri nets (S.~Gaubert and J.~Mairesse~\cite{GaubertMairesseIEEE}) and some queuing networks (J. Mairesse~\cite{Mairesse}, B.~Heidergott~\cite{CaractMpQueuNet}) and many concrete applications. Let us cite job-shops models (G.~Cohen~et~al.\cite{cohen85a}), train networks (H.~Braker~\cite{Braker}, A.~de Kort and B.~Heidergott~\cite{RailwayMpKHA}), computer networks (F.~Baccelli~\cite{TCPmp}) or a statistical mechanics model (R.~Griffiths~\cite{Griffiths}). For more details about modelling, the reader is referred to the books by F.~Baccelli and al.~\cite{BCOQ} and by B.~Heidergott and al.~\cite{MpAtWork}.

\subsection{Law of large numbers}
The sequences satisfying Equation~(\ref{eqdefx}) have been studied in many papers. Law of large numbers have been proved among others by J.E.~Cohen~\cite{Cohen}, F.~Baccelli and Z.~Liu~\cite{BaccelliLiu}, and more recently by T.~Bousch and J.~Mairesse~\cite{BouschMairesse}. 

If matrix $A$ has no line of $-\infty$, then $x\mapsto Ax$ is 1-Lipschitz for the supremum norm. Therefore, we can assume $x_0=0$, and we do it from now on.

T.~Bousch and J.~Mairesse have proved (cf.~\cite{BouschMairesse}) that, if $A(0)0$ is integrable, then the sequence $\left(\frac{1}{n}y(n,0)\right)_{n\in\N}$ converges almost-surely and in mean. The proof is still true if $\max_{ij}A_{ij}^+(0)$ is integrable and the limit can be $-\infty$.

Therefore, the sequence $\left(\frac{1}{n}x(n,0)\right)_{n\in\N}$ converges in law. But, it does not necessary converges almost-surely, as illustrated by examples bellow. T.~Bousch and J.~Mairesse have also proved that, if $\max_{A_{ij}(0)\neq -\infty}|A_{ij}(0)|$ is integrable, then $\left(\frac{1}{n}x(n,0)\right)_{n\in\N}$ convergences almost-surely if and only if the limit of $\left(\frac{1}{n}y(n,0)\right)_{n\in\N}$ is constant.  Thanks to their former results, we give necessary conditions for $\left(\frac{1}{n}x(n,x_0)\right)_{n\in\N}$ to converge almost surely under the usual integrability condition.

F.~Baccelli has proved (cf.~\cite{Baccelli}) by induction on the size of the matrices, that $\left(\frac{1}{n}x(n,0)\right)_{n\in\N}$ converges almost-surely under the following additional hypotheses: each entry of the matrix is either almost-surely  $-\infty$, or almost-surely non-negative and the diagonal entries are non-negative (precedence condition).  T.~Bousch and J.~Mairesse have proved (cf.~\cite{BouschMairesse}) that, if $\max_{A_{ij}(0)\neq -\infty}|A_{ij}(0)|$ is integrable, then the precedence condition is sufficient. Practically, the limit of $\left(\frac{1}{n}y(n,0)\right)_{n\in\N}$ is constant if the diagonal coefficients are almost-surely finite.  We use the induction method without the additional hypotheses. 
We obtain a general scheme (Theorem~\ref{thgene}) to prove that $\left(\frac{1}{n}x(n,x_0)\right)_{n\in\N}$ converges almost surely.
When $\max_{A_{ij}(0)\neq -\infty}|A_{ij}(0)|$ is integrable, it allows us:
\begin{enumerate}[-]
\item to give a necessary and sufficient condition for the convergence of $\left(\frac{1}{n}x(n,0)\right)_{n\in\N}$ when the sequence $\sAn$ is i.i.d.,
\item to prove that, if $\sAn$ satisfy a condition of reinforced ergodicity and a condition of fixed structure (i.e. $\P\left(A_{ij}(0)=-\infty\right)\in\{0,1\}$), then  $\left(\frac{1}{n}x(n,0)\right)_{n\in\N}$ converges almost-surely,
\item and to reprove the convergence of $\left(\frac{1}{n}x(n,0)\right)_{n\in\N}$ if the diagonal entries are never~$-\infty$.
\end{enumerate}

In the next section, we will state our results, then we will prove a necessary condition for the convergence of $\left(\frac{1}{n}x(n,0)\right)_{n\in\N}$ (Theorem~\ref{thCN}), then the sufficient condition (Theorem~\ref{thgene}) and finally its three consequences.
\section{Presentation of the results}
The first result is the following, which directly follows from Kingman's theorem and can be traced back to J.E.~Cohen~\cite{Cohen}:
\begin{th-def}[Maximal Lyapunov exponent]\label{defexptop}\ \\
If $\sAn$ is an ergodic sequence of random matrices with entries in $\rma$, such that $\max_{ij} A_{ij}^+(0)$ is integrable, then the sequences $\left(\frac{1}{n} \max_i x_i(n,0)\right)_{n\in\N}$ and $\left(\frac{1}{n} \max_i y_i(n,0)\right)_{n\in\N}$ converges almost-surely to the same constant $\gamma\in\rma$, which is called maximal Lyapunov exponent of $\left(A(n)\right)_{n\in\N}$.

We denote this constant by $\gamma\left(\left(A(n)\right)_{n\in\N}\right)$, or  $\gamma(A)$.
\end{th-def}
\begin{rem}
The constant $\gamma(A)$ is well-defined even if $\sAn$ has a line of~$-\infty$.\\
The variable $\max_i x_i(n,0)$ is equal to $\max_{ij} A^n_{ij}$.
\end{rem}
Let us associate to our sequence of random matrices a graph to split the problem. We also set the notations for the rest of the text.
\begin{defi}[Graph of possible incidences]\label{defGiA}
For every $x\in\rma^{[1,\cdots,d]}$ and every subset $I\subset [1,\cdots,d]$, we write:
$$x^I:=(x_i)_{i\in I}.$$
Let $\sAn$  be a stationary sequence of random matrices with values in $\rma^{d\times d}$. 
\begin{enumerate}[i)]
\item  The graph of possible incidences of $\sAn$, denoted by $\G(A)$, is the directed graph whose nodes are the integers between 1 and d and whose arcs are the pairs $(i,j)$ such that $\P(A_{ij}(0)\neq-\infty)>0$.

\item We denote by $c_1,\cdots , c_K$ the strongly connected components of $\G(A)$. In the sequel, we just say components of $\G(A)$.

To each component $c_m$, we associate the following elements:
$$A^{(m)}(n):=(A_{ij}(n))_{i,j\in c_m}\textrm{, }\gamma^{(m)}:=\gamma(A^{(m)})\textrm{, }$$
$$x^{(m)}(n,x_0):=(A^{(m)})^n(x_0)^{c_m}
 \textrm{ and } y^{(m)}(n,x_0):=(A^{(m)})^{-n}(x_0)^{c_m} $$

\item A component $c_l$ is reachable from a component $c_m$, if $m=l$ or if there exists a path on $\G(A)$ from a node in $c_m$ to a node in  $c_l$. In this case, we write $m\rightarrow l$.

To each component $c_m$, we associate the following elements:
$$E_m:=\{l\in [1,\cdots,K]|m\rightarrow l\}
\textrm{, }
\gamma^{[m]}:=\max_{l\in E_m}\gamma^{(l)}
\textrm{, }$$
$$F_m:=\bigcup_{l\in E_m}c_l 
\textrm{, }
A^{[m]}(n):=\left(A_{ij}(n)\right)_{i,j\in F_m}$$
$$x^{[m]}(n,x_0):=\left(A^{[m]}\right)^n(x_0)^{F_m} \textrm{ and }
y^{[m]}(n,x_0):=\left(A^{[m]}\right)^{-n}(x_0)^{F_m}.$$

\item A component $c_m$ is final (or source, in the terminology of discrete event systems) if $E_m=\{m\}$, that is if, for every $l\in [1,\cdots,K]$, we have:
$$m\rightarrow l\Rightarrow l=m.$$
It is initial if, for every $l\in [1,\cdots,K]$, we have:
$$l\rightarrow m\Rightarrow l=m.$$
A component is said to be trivial, if it has only one node $i$ and $\P(A_{ii}(1)\neq -\infty)=0$.
\item
To each component $c_m$, we associate the following sets: 
$$G_m:=\{l\in E_m|\exists p\in [1,\cdots,K], m\rightarrow l\rightarrow p, \gamma^{(p)}=\gamma^{[m]}\},$$
$$H_m:=\bigcup_{l\in G_m}c_l\textrm{ , }  A^{\{m\}}(n):=\left(A_{ij}(n)\right)_{i,j\in H_l} $$
$$x^{\{m\}}(n,x_0):=\left(A^{\{m\}}\right)^{n}(x_0)^{H_m}\textrm{ and } y^{\{m\}}(n,x_0):=\left(A^{\{m\}}\right)^{-n}(x_0)^{H_m}.$$
\item A component $c_m$ is called dominating if $G_m=\{m\}$, that is if for every $l\in E_m\backslash\{m\}$, we have: 
$\gamma^{(m)}>\gamma^{(l)}.$
\end{enumerate}
\end{defi}
\begin{rem}[Paths on $\G(A)$]
Equation~(\ref{eqdefpdt}) can be read as '$A^{n}_{ij}$ is the maximum of the weights of paths from $i$ to $j$ with length $n$ on $\G(A)$, the weight of the $k^\textrm{th}$ arc being given by $A(-k)$'. Thus $y_i(n,0)$  is the maximum of the weights of paths on $\G(A)$ with initial node~$i$ and length $n$. The coefficients $y^{(m)}_i(n,0)$ and $y^{\{m\}}_i(n,0)$ are the maximum of the weights of paths on the subgraph of $\G(A)$ with nodes in $c_m$ and $H_m$ respectively.

Consequently $\gamma^{(m)}$ is the average maximal weight of path on $c_m$
\end{rem}\ \\

The first new result is a necessary condition for $x(n,X_0)$ to satisfy a strong law of large numbers:
\begin{theo}\label{thCN}
Let $\left(A(n)\right)_{n\in\N}$ be a stationary and ergodic sequence of random matrices  with values in $\rma^{d\times d}$ and no line of~$-\infty$, such that $\max_{ij} A_{ij}^+(0)$ is integrable.

If the limit of $\left(\frac{1}{n}y(n,0)\right)_{n\in\N}$ is deterministic, then it is given by:
\begin{equation}\label{eqcvgcecompy}
\forall m\in [1,K], \lim_n\frac{1}{n}y^{c_m}(n,0)=\gamma^{[m]}\1~\mathrm{ a.s.},
\end{equation}
where $\1$ is the vector whose coordinates are all~$1$.\\
That being the case, for every component $c_m$ of $\G(A)$, $A^{\{m\}}(0)$ has no line of~$-\infty$.\\

If $\left(\frac{1}{n}x(n,0)\right)_{n\in\N}$ converges almost-surely, then its limit is deterministic and is equal to that of $\left(\frac{1}{n}y(n,0)\right)_{n\in\N}$, that is we have:
\begin{equation}\label{eqcvgcecompx}
\forall m\in [1,K], \lim_n\frac{1}{n}x^{c_m}(n,0)=\gamma^{[m]}\1~\mathrm{ a.s.},
\end{equation}
\end{theo}

The following theorem gives a scheme to prove converse theorems:
\begin{theo}\label{thgene}
Let $\left(A(n)\right)_{n\in\N}$ be an ergodic sequence of random matrices with values in $\rma^{d\times d}$ that satisfy the three following hypotheses:
\begin{enumerate}
\item for every component $c_m$ of $\G(A)$, $A^{\{m\}}(0)$ has no line of~$-\infty$.
\item for every dominating component $c_m$ of $\G(A)$, $\lim_n\frac{1}{n}y^{(m)}(n,0)=\gamma^{(m)}\1~\mathrm{ a.s.\,.}$
\item for every subsets $I$ and $J$  of $[1,\cdots,d]$, such that random matrices $\tilde{A}(n)=\left(A_{ij}(n)\right)_{i,j\in I\cup J}$ has no line of~$-\infty$ and split along $I$ and $J$ following the equation
\begin{equation}\label{decompblocs}
\tilde{A}(n)=:\left(\begin{array}{cc}
B(n)&D(n)\\
-\infty&C(n)
\end{array}\right),
\end{equation}
such that $\G(B)$ is strongly connected and  $D(n)$ is not almost surely $(-\infty)^{I\times J}$, we have:
\begin{equation}\label{eqexistchem}
\P\left(\left\{\exists i\in I, \forall n\in\N, \left(B(-1)\cdots B(-n)D(-n-1)0\right)_i=-\infty\right\}\right)=0.
\end{equation}
\end{enumerate}
Then the limit of $\left(\frac{1}{n}y(n,0)\right)_{n\in\N}$ is given by Equation~(\ref{eqcvgcecompy}).\\

If Hypothesis~1. is strengthened by demanding that $A^{\{m\}}(0)0$ is integrable, then the sequence $\left(\frac{1}{n}x(n,0)\right)_{n\in\N}$ converges almost-surely and its limit is given by Equation~(\ref{eqcvgcecompx}).
\end{theo}

\begin{rem}[Paths on $\G(A)$, continued]
Let us interpret the three hypotheses with the paths on~$\G(A)$.
\begin{enumerate}
\item The hypothesis on $A^{\{m\}}(0)$ means that, whatever the initial condition~\mbox{$i\in c_m$,} there is always a path beginning in~$i$ and staying in~$H_m$.
\item The hypothesis on dominating component means that, whatever the initial condition~$i$ in dominating component $c_m$, there is always a path beginning in~$i$ with average weight $\gamma^{(m)}$. It is necessary, as can be shown by an method analogeous to that of~\cite{Baccelli}.
\item We will use the last hypothesis with $\tilde{A}(n)=A^{\{m\}}(n)$, $B(n)=A^{(m)}(n)$.  It means there is a path from~$i\in c_m$, to $H_m\backslash c_m$.
Once we know that the limit of $\left(\frac{1}{n}y(n,0)\right)_{n\in\N}$ is given by Equation~(\ref{eqcvgcecompy}) this hypothesis is obviously necessary when $\gamma^{(m)}< \gamma^{[m]}$.
\end{enumerate}

Thanks to hypotheses~1. and~3., for every component~$c_m$ and every node~$i\in c_m$, there is always a path beginning in~$i$ reaching a dominating component $c_k$ with Lyapunov exponent~$\gamma^{(k)}=\gamma^{[m]}$  and staying in that component. Thanks to Hypothesis~2., this paths has $\gamma^{[m]}$ as average weight.

It remains to prove that there is no path in $F_m$ with average weight strictly greater than~$\gamma^{[m]}$ and goes from Equation~(\ref{eqcvgcecompy}) to Equation~(\ref{eqcvgcecompx}). This is possible thanks to theorem~\ref{decomplyap}, from~\cite{Baccelli} and theorem~\ref{thvincent}, from~\cite{vincent} respectively.
\end{rem}

The three announced results follow from this scheme:
\begin{theo}[Independent case]\label{thiid1}
If $\left(A(n)\right)_{n\in\N}$ is a sequence of i.i.d. random matrices with values in $\rma^{d\times d}$ and no line of~$-\infty$, such that $\max_{A_{ij}(0)\neq -\infty}|A_{ij}(0)|$ is integrable, then the sequence $\left(\frac{1}{n}x(n,0)\right)$ converges almost-surely\footnote{Actually, the convergence of $\left(\frac{1}{n}x(n,0)\right)$ is always proved under the -- slightly weaker but much more difficult to check -- condition $\forall m, A^{\{m\}}(0)0\in\L^1$, which appears in Theorem~\ref{thgene}.\label{noteinteg}}
if and only if for every component $c_m$, $A^{\{m\}}$ has no line of~$-\infty$.
That being the case the limit is given by Equation~(\ref{eqcvgcecompx}).
\end{theo}

\begin{theo}[Fixed structure case]\label{thSF1}
If $(\Omega,\theta,\P)$ is a measurable dynamical system, and $A:\Omega\rightarrow\rma^{d\times d}$ is a random matrix with no line of~$-\infty$, such that:
\begin{enumerate}
\item for every $i,j$, $A_{ij}$ is integrable or almost-surely $-\infty$,
\item for every $k\le d$, $\theta^k$ is ergodic,
\end{enumerate}
then, the sequence $\left(\frac{1}{n}x(n,0)\right)_{n\in\N}$ associated to the sequence $\sAn=\left(A\circ\theta^n\right)_{n\in\N}$ converges almost-surely\hbox{$^{\ref{noteinteg}}$} and its limit is given by Equation~(\ref{eqcvgcecompx}).
\end{theo}

\begin{theo}[Precedence case]\label{thprec}
If $\left(A(n)\right)_{n\in\N}$ is a stationary and ergodic sequence of random matrices  with values in $\rma^{d\times d}$ such that $\max_{A_{ij}(0)\neq -\infty}A^+_{ij}(0)$ is integrable and for every $i\in [1,\cdots,d]$
\begin{equation}\label{eqpreced}
\P\left(A_{ii}(0)=-\infty\right)=0,
\end{equation}
then the limit of $\left(\frac{1}{n}y(n,0)\right)$ is given by Equation~(\ref{eqcvgcecompy}).

If condition~(\ref{eqpreced}) is strengthened by demanding that $A_{ii}(0)$ is integrable, then the sequence $\left(\frac{1}{n}x(n,0)\right)$ converges almost-surely\hbox{$^{\ref{noteinteg}}$} and its limit is given by Equation~(\ref{eqcvgcecompx}).
\end{theo}

To end this section, we give two examples that show that neither the fixed structure, neither the independence ensure the strong law of large numbers.
\begin{ex}[\cite{theseMairesse}]\label{exechanges}
Let us set $ \Omega=\{ \omega_0, \omega_1\}$ and $\P=\frac{1}{2}(\delta_{\omega_0}+\delta_{\omega_1})$.
Let  $\theta$ be the function that exchange $\omega_0$ and $\omega_1$ and  $A$, from $ \Omega$ to $\rma^{2\times 2}$ be defined by 
$$A( \omega_0)=\left(\begin{array}{cc}-\infty&0\\ 0&-\infty\end{array}\right) \textrm{ and } A( \omega_1)=\left(\begin{array}{cc}-\infty&1\\ 0&-\infty\end{array}\right).$$
It can be checked that $( \Omega,\theta,\P)$ is an ergodic dynamical system, and that the sequence $\left(A(n)\right)_{n \in\N}$ has fixed structure. Moreover $\G(A)$ is strongly connected.
The sequence is a degenerate Markov chain: $A\circ\theta^n=A( \omega_0)\Leftrightarrow A\circ\theta^{n+1}=A( \omega_1)$. 

The multiplication by $A( \omega_0)$ exchange the coordinates, and the multiplication by  $A( \omega_1)$ does the same, and then increases the first coordinate by~$1$. Therefore the sequence is defined by the following equations:
\begin{equation}\label{eqexechanges}
\begin{array}{rclcrcl}
x(2n,z)( \omega_0)&=&(z_1+n,z_2)'&\textrm{et}&x(2n+1,z)( \omega_0)&=&(z_2,z_1+n)'\\
x(2n,z)( \omega_1)&=&(z_1,z_2+n)'&\textrm{et}&x(2n+1,z)( \omega_0)&=&(z_2+n+1,z_1)'.\\
\end{array}
\end{equation}
Therefore the sequence $\left(\frac{1}{n}x(n,0)\right)_{n\in\N}$ almost-surely does not converges.\\

As $A(n)$ has fixed structure and $\G(A)$ has only one component, it proves that Hypothesis~2. of Theorem~\ref{thiid1} is necessary.
\end{ex}

\begin{ex}[\cite{BouschMairesse2}]\label{exmairesse}
Let $\left(A(n)\right)_{n \in\N}$ be the  sequence of i.i.d. random variables taking values
$$B=\left(\begin{array}{ccc}0&-\infty&-\infty\\ 0&-\infty&-\infty\\ 0&1&1\end{array}\right) 
\textrm{ and } 
C=\left(\begin{array}{ccc}0&-\infty&-\infty\\ 0&-\infty&0\\ 0&0&-\infty\end{array}\right)$$
with probabilities $p>0$ and $1-p>0$.
Let us compute the action of $B$ and $C$ on vectors of type $(0,x,y)'$, with $x,y\ge0$:
$$B(0,x,y)'=(0,0,\max(x,y)+1)'\textrm{ and } C (0,x,y)'=(0,y,x)'.$$
Therefore $x_1(n,0)=0$ and $\max_ix_i(n+1,0)=\#\{0\le k\le n|A(k)=B\}$. Practically, if 
$A(n)=B$, then $x(n+1,0)=\left(0,0,\#\{0\le k\le n|A(k)=B\}\right)'$, and if $A(n)=C$ and $A(n-1)=B$, then $x(n+1,0)=\left(0,\#\{0\le k\le n|A(k)=B\},0\right)'$. Since $\left(\frac{1}{n}\#\{0\le k\le n|A(k)=B\}\right)_{n\in\N}$ converges almost-surely to $p$, we see:
\begin{equation}\label{eqexmairesse}
\begin{array}{c}
\lim_n \frac{1}{n} x_1(n,0)=0~\mathrm{ a.s.}\\
\forall i\in\{2,3\}, \liminf_n \frac{1}{n} x_i(n,0)=0\textrm{ and }\limsup_n \frac{1}{n} x_i(n,0)=p\textrm{  a.s.\,.}
\end{array}
\end{equation}
Therefore the sequence $\left(\frac{1}{n}x(n,0)\right)_{n\in\N}$ almost-surely does not converge.\\

We notice that $\G(A)$ has two components $c_1=\{1\}$ and $c_2=\{2,3\}$, with Lyapunov exponents $\gamma^{(1)}=0$ and $\gamma^{(2)}=p$, and $2\rightarrow 1$. Therefore we check that $A^{\{2\}}(n)$ has a line of $-\infty$ with probability $p$.
\end{ex}\ \\

The last example shows the necessity of the integrability conditions in the former theorems: it satisfy every hypothesis of each the three theorems, except for the integrability conditions, but the associated $\left(x(n,0)\right)_{n\in\N}$ does not satisfy a strong law of large numbers.
\begin{ex}[Integrability]\label{exinteg}
Let $\left(X_n\right)_{n\in\N}$ be a sequence of real variables satisfying $X_n\ge 1\textrm{ a.s.}$ and $\E(X_n)=+\infty$. The sequence of matrices is defined by:
$$A(n)=\left(\begin{array}{ccc}
-X_n&-X_n&0\\
-\infty&0&0\\
-\infty&-\infty&0
\end{array}\right)$$
A straightforward computation shows that $x(n,0)=\left(\max(-X_n,-n),0,-n\right)'$ and $y(n,0)=\left(\max(-X_0,-n),0,-n\right)'$. Since $\P\left(\lim_n\frac{1}{n}X_n=0\right)=0$, it implies that $\left(\frac{1}{n}x(n,0)\right)_{n\in\N}$ converges to $(0,0,-1)'$ in probability but not almost-surely.\\

Let us notice that the limit of $\left(\frac{1}{n}y(n,0)\right)_{n\in\N}$ can also be found by applying one of the theorems  and computing that each component has exactly one node and $\gamma^{(1)}=-\E(X_n)=-\infty$, $\gamma^{(2)}=0$ and $\gamma^{(3)}=-1$.
\end{ex}

\section{Proofs}
\subsection{Necessary conditions}

\subsubsection{Formula for the limit}
Let us denote by $L$ the limit of $\left(\frac{1}{n}y(n,0)\right)_{n\in\N}$, which exists, according to~\cite{BouschMairesse}, and is assumed to be deterministic.

By definition of $\G(A)$, if $(i,j)$ is an arc of $\G(A)$, then, with positive probability, we have:
 $$L_i=\lim_n \frac{1}{n}y_i(n,0)\ge \lim_n \frac{1}{n} (A_{ij}(-1)+y_j(n,0)\circ\theta^{-1})=0+ L_j\circ\theta^{-1}=L_j.$$

If $m \rightarrow p$, then for every $i\in c_m$ and $j\in c_p$, there exists a path on $\G(A)$ from $i$ to $j$, therefore $L_i\ge L_j$. Since this is true for every $j\in F_m$, we have:
\begin{equation}\label{eqlim}
L_i =\max_{j\in F_m} L_j
\end{equation}
To show that $\max_{j\in F_m} L_j=\gamma^{[m]}$, we have to study the Lyapunov exponents of sub-matrices.
The following proposition states some easy consequences of definition~\ref{defGiA}, which will be useful in the sequel.
\begin{prop}\label{propGiA}
The notations are those of definition~\ref{defGiA}\,.
\begin{enumerate}[i)]
\item for every $m\in [1,\cdots,K]$, $x^{[m]}(n,x_0)=x^{F_m}(n,x_0)$.
\item for every  $m\in [1,\cdots,K]$,  and every $i\in c_m$, we have:
$$x_i^{c_m}(n,0)= x_i^{[m]}(n,0)\ge x_i^{\{m\}}(n,0)\ge x_i^{(m)}(n,0).$$
\begin{equation}\label{OrdreDesy}
y_i^{c_m}(n,0)= y_i^{[m]}(n,0)\ge y_i^{\{m\}}(n,0)\ge y_i^{(m)}(n,0).
\end{equation}
\junk{$$\max_{i_l\in [1,d]}\sum_{l=0}^{n-1}A_{i_{l+1}i_l}(l) =\max_{i_l\in F_m}\sum_{l=0}^{n-1}A_{i_{l+1}i_l}(l) \ge \max_{i_l\in H_m}\sum_{l=0}^{n-1}A_{i_{l+1}i_l}(l) \ge \max_{i_l\in c_m}\sum_{l=0}^{n-1}A_{i_{l+1}i_l}(l).$$}
\item relation $\rightarrow$ is a partial order. Initial and final components are  minimal and maximal elements for this order.
\item If $A(0)$ has no line of~$-\infty$, then for every $m\in [1,\cdots,K]$, $A^{[m]}(0)$ has no line of~$-\infty$.
Practically final components has no line of~$-\infty$ and are never trivial.
\item for every $l\in E_m$, we have $\gamma^{(l)}\le\gamma^{[l]}\le\gamma^{[m]}$ and \mbox{$G_m=\{l\in E_m|\gamma^{[l]}=\gamma^{[m]}\}.$}
\end{enumerate} 
\end{prop}

The next result is about Lyapunov exponents. It is already in~\cite{BaccelliLiu} and its proof does not uses the additional hypotheses of this article. For a point by point checking, the reader is referred to~\cite{theseGM}.
\begin{theo}[\cite{Baccelli}]\label{decomplyap}
If $\sAn$ is a stationary and ergodic sequence of random matrices with values in $\rma^{d\times d}$ such that $\max_{i,j} A_{ij}^+$ is integrable, then $\gamma(A)=\max_l\gamma^{(l)}$.
\end{theo}

Applying this theorem to sequences $\left(A^{[m]}(n)\right)_{n\in\N}$ and $\left(A^{\{m\}}(n)\right)_{n\in\N}$, we obtain the following proposition
\begin{prop}\label{propGiA2}
For every $m\in [1,\cdots,K]$, we have $\gamma(A^{\{m\}})=\gamma(A^{[m]})=\gamma^{[m]}$.
\end{prop}

It follows from proposition~\ref{propGiA} and  the definition of Lyapunov exponents  that for every component $c_m$ of $\G(A)$,
$$\max_{i\in F_m}L_i=\lim_n\frac{1}{n}\max_{i\in F_m}y_i(n,0)=\gamma(A^{\{m\}}).$$

Combining this with Equation~(\ref{eqlim}) and proposition~\ref{propGiA2}, we deduce that the limit of $\left(\frac{1}{n} y(n,0)\right)_{n\in\N}$ is given by Equation~(\ref{eqcvgcecompy}).

\subsubsection{\texorpdfstring{$A^{\{m\}}(0)$}{$A^{\{m\}}(0)$} has no line of $-\infty$}
We still have to show that for every component $c_m$, $A^{\{m\}}(0)$ has no line of~$-\infty$. Let us assume it has one. Therefore, there exists $m\in [1,\cdots,d]$ and $i\in c_m$ such that the set $$\{\forall j\in H_m, A_{ij}(-1)=-\infty\}$$ has positive probability. On this set, we have:
$$y_i(n,0)\le \max_{j\in F_m\backslash H_m} A_{ij}(-1) +\max_{j\in F_m\backslash H_m} y_j(n-1,0)\circ\theta^{-1}.$$
Dividing by $n$ and letting $n$ to $+\infty$, we have $L_i\le\max_{j\in F_m\backslash H_m} L_j$, which, because of Equation~(\ref{eqcvgcecompy}) becomes $\gamma^{[m]}\le\max_{k\in E_m\backslash G_m} \gamma^{[k]}$. This last inequality contradicts proposition~\ref{propGiA}~$v)$. Therefore the hypothesis that $A^{\{m\}}(0)$ has a line of~$-\infty$.

\subsubsection{The limit is constant}
Let us assume that $\left(\frac{1}{n}x(n,0)\right)_{n\in\N}$  converges almost-surely to a limit $L'$. Up to a change of probability space, we can assume that $A(n)=A\circ\theta^n$, where $A$ is a random variable  and $(\Omega,\theta,\P)$ is an invertible ergodic  measurable dynamical system.

It follows from~\cite{BouschMairesse} that $\left(\frac{1}{n}y(n,0)\right)_{n\in\N}$ converges almost-surely and  
$$\frac{1}{n}y(n,0)-\frac{1}{n+1}y(n+1,0)\stackrel{\P}{\rightarrow} 0.$$
We compound each term of this relation by $\theta^{n+1}$ and, since $x(n,0)=y(n,0)\circ\theta^n $, it proves that:
$$\frac{1}{n}x(n,0)\circ\theta-\frac{1}{n+1}x(n+1,0)\stackrel{\P}{\rightarrow}0.$$
When $n$ tends to $+\infty$, it  becomes $L'\circ\theta -L'=0$. Since $\theta$ is ergodic, this implies that $L'$ is constant.\\

 Since $\frac{1}{n}y(n,0)=\frac{1}{n}x(n,0)\circ\theta^n$, $L'$ and $L$ have the same law. Since $L'$ is constant, $L=L'$ almost-surely, therefore $L$ is also the limit of $\left(\frac{1}{n}x(n,0)\right)_{n\in\N}$.
This proves formula~(\ref{eqcvgcecompx}) and concludes the proof of theorem~\ref{thCN}\,.

\subsection{Main theorem}
\subsubsection{Right products}
In this section, we prove Theorem~\ref{thgene}. We begin with the result on $y(n,0)$.

It follows from propositions~\ref{propGiA} and~\ref{propGiA2} and the definition of Lyapunov exponents that we have, for every component $c_m$ of $\G(A)$,
\begin{equation}\label{eqmajor}
\limsup_n\frac{1}{n}y^{c_m}(n,0)\le\gamma^{[m]}\1~\mathrm{a.s.\,.}
\end{equation}

Therefore, it is sufficient to show that $\liminf_n\frac{1}{n}y^{c_m}(n,0)\ge\gamma^{[m]}\1~\mathrm{a.s.}$\,. Because of proposition~\ref{propGiA}~$i)$, it is sufficient to show that 
\begin{equation}\label{eqrec}
\lim_n\frac{1}{n}y^{\{m\}}(n,0)=\gamma^{[m]}\1.
\end{equation}
We prove Equation~(\ref{eqrec}) by induction on the size of $G_m$. The initialization of the induction is exactly Hypothesis~$2.$ of Theorem~\ref{thgene}.

Let us assume that Equation~(\ref{eqrec}) is satisfied by every $m$ such that the size of $G_m$ is less than $N$, and let $m$ be such that the size of $G_m$ is~$N+1$. Let us take $I=c_m$ and $J=H_m\backslash c_m$. If $c_m$ is not trivial, it is  the situation of Hypothesis~$3.$ with $\tilde{A}=A^{\{m\}}$, which has no line of~$-\infty$ thanks to Hypothesis~$1.$\,. Therefore  Equation~(\ref{eqexistchem}) is satisfied. If $c_m$ is trivial, $\G(B)$ is not strongly connected, but Equation~(\ref{eqexistchem}) is still satisfied because $D(-1)0=(\tilde{A}(-1)0)^I\in\R^I$.

Moreover $J$ is the union of the $c_k$ such that $k\in G_m\backslash\{m\}$, thus the  induction hypothesis implies that:
$$\forall j\in J, j\in c_k\Rightarrow\lim_n\frac{1}{n}(C^{-n} 0)_j=\lim_n\frac{1}{n}y^{\{k\}}_j(n,0)=\gamma^{[k]}~\mathrm{ a.s. }. $$
Because of proposition~\ref{propGiA2}~$ii)$, $\gamma^{[k]}=\gamma^{[m]},$ therefore the right side of the last equation is $\gamma^{[m]}$ and we have:
\begin{equation}\label{eqrecJ}
\lim_n\frac{1}{n}(y^{\{m\}})^J(n,0)=\lim_n \frac{1}{n}C^{-n}0= \gamma^{[m]}\1~\mathrm{ a.s. }.
\end{equation}

Now Equation~(\ref{eqexistchem}) ensures that for every $i\in I$, there exists almost-surely a $T\in \N$ and a $j\in J$ such that $\left(B(-1)\cdots B(-T)D(-T-1)\right)_{ij}\neq-\infty$. Since we have $\lim_n\frac{1}{n}\left(C(-T)\cdots C(-n)0\right)_{j}=\gamma^{[m]} \textrm{ a.s.}$, it implies that:
\begin{eqnarray*}
\lefteqn{\liminf_n\frac{1}{n} y_i^{\{m\}}(n,0)}\\
&\ge& \lim_n \frac{1}{n}\left(B(-1)\cdots B(-T)D(-T-1)\right)_{ij}+ \lim_n\frac{1}{n}\left(C(-T)\cdots C(-n)0\right)_{j}=\gamma^{[m]}~\mathrm{ a.s. }
\end{eqnarray*}
Because of upper bound~(\ref{eqmajor}) and inequality~(\ref{OrdreDesy}), it implies that
$$\lim_n\frac{1}{n}(y^{\{m\}})^I(n,0)= \gamma^{[m]}\1~\mathrm{ a.s. }.,$$
which, because of Equation~(\ref{eqrecJ}), proves Equation~(\ref{eqrec}). This  concludes the induction and  the proof of the result on~$y(n,0)$.

\subsubsection{Left products}
To deduce the results on $x(n,0)$ from those on $y(n,0)$, we introduce the following theorem-definition, which is a special case of the main theorem of J.~M.~Vincent~\cite{vincent} and  directly follows from Kingman's theorem:

\begin{th-def}[\cite{vincent}]\label{thvincent}
If $\left(A(n)\right)_{n\in\N}$ is a stationary and ergodic sequence of random matrices  with values in $\rma^{d\times d}$ such that $A(0)0$ is integrable, then there are two real numbers $\gamma(A)$ and $\gamma_b(A)$ such that 
$$\lim_n\frac{1}{n}\max_ix_i(n,0)=\frac{1}{n}\max_iy_i(n,0)=\gamma(A) ~\mathrm{ a.s.} $$
$$\lim_n\frac{1}{n}\min_ix_i(n,0)=\frac{1}{n}\min_iy_i(n,0)=\gamma_b(A) ~\mathrm{ a.s.} $$
\end{th-def}
It implies the following corollary, which makes the link between the results on~$y(n,0)$ and those on~$x(n,0)$ when all $\gamma^{[m]}$ are equal, that is when $\gamma(A)=\gamma_b(A)$.
\begin{coro}\label{lempassgauche}
If $\left(A(n)\right)_{n\in\N}$ is a stationary and ergodic sequence of random matrices  with values in $\rma^{d\times d}$ such that $A(0)0$ is integrable then $$ \lim_n\frac{1}{n}x(n,0)=\gamma(A)\1\textrm{ if and only if }\lim_n\frac{1}{n}y(n,0)=\gamma(A)\1.$$
\end{coro}

Let us go back to the proof of the general result on $x(n,0)$.
Because of propositions~\ref{propGiA} and~\ref{propGiA2} and  the definition of Lyapunov exponents, we already have, for every component $c_m$ of $\G(A)$,
$$\limsup_n\frac{1}{n}x^{c_m}(n,0)\le\gamma^{[m]}\1~\mathrm{a.s.\,.}$$

Therefore it is sufficient to show that $\liminf_n\frac{1}{n}x^{c_m}(n,0)\ge\gamma^{[m]}\1~\mathrm{a.s.\,.}$ and even that
$$\lim_n\frac{1}{n}x^{\{m\}}(n,0)=\gamma^{[m]}\1.$$
Because of corollary~\ref{lempassgauche}, it is equivalent to $\lim_n\frac{1}{n}y^{\{m\}}(n,0)=\gamma^{[m]}\1.$
Since all components of $\G(A^{\{m\}})$ are components of $\G(A)$ and have the same Lyapunov exponent~$\gamma^{[m]}$, it follows from the result on the~$y(n,0)$ applied to $A^{\{m\}}$.

\subsection{Independent case}
In this section, we prove Theorem~\ref{thiid1}. 

Because of Theorem~\ref{thCN}, it is sufficient to show that, if $\left(A(n)\right)_{n\in\N}$ is a sequence of i.i.d. random matrices  with values in $\rma^{d\times d}$ such that $\max_{A_{ij}(0)\neq -\infty}|A_{ij}(0)|$ is integrable and for every component $c_m$, $A^{\{m\}}$ has no line of~$-\infty$, then the sequence $\left(\frac{1}{n}x(n,0)\right)$ converges almost-surely. To do this, we will prove that in this situation, the hypotheses of Theorem~\ref{thgene} are satisfied. Hypothesis~$1.$ is exactly Hypothesis~$1.$ of Theorem~\ref{thiid1} and hypotheses $2.$ and $3.$ readily follow from the next theorem and lemma respectively.
\begin{theo}[D.~Hong~\cite{Hong}]\label{thHong}
If $\left(A(n)\right)_{n\in\N}$ is a sequence of i.i.d. random matrices with values in $\rma^{d\times d}$ such that $A(1)0$ is integrable, $A(1)$ has no line of~$-\infty$ and $\G(A)$ is strongly connected, then $\gamma(A)=\gamma_b(A)$.
\end{theo}

\begin{lem}\label{lemMifini}\ 
Let $\sAn$ be a stationary sequence of random matrices with values in $\rma^{d\times d}$ with no line of~$-\infty$. Let us assume that there exists a partition $(I,J)$ of $[1,\cdots,d]$ such that $A=\tilde{A}$ satisfy Equation~(\ref{decompblocs}), with $\G(B)$ strongly connected.
For every $i\in I$, let us define
$$\A_i:=\left\{\forall n\in\N, \left(B(1)\cdots B(n)D(n+1)0\right)_i=-\infty\right\}.$$

\begin{enumerate}
\item If $\omega\in\A_i$, then we have $\forall n\in\N, \exists i_n\in I \left(B(1)\cdots B(n)\right)_{ii_n}\neq -\infty.$
\item If the random matrices $A(n)$ are i.i.d., and if $\P\left(D=(-\infty)^{I\times J}\right)<1$, then for every $i\in I$, we have $ \P(\A_i)=0.$
\end{enumerate}
\end{lem}

\begin{proof}\ 
\begin{enumerate}
\item For every $\omega\in\A_i$, we prove our result by induction on $n$.

Since the $A(n)$ have no line of~$-\infty$, there exists an $i_1\in [1,\cdots,d]$, such that $A_{ii_1}(1)\neq -\infty$.
Since $\left(D(1)0\right)_{i}=-\infty$, every entry on line~$i$ of~$D(1)$ is $-\infty$, that is $A_{ij}(1)=-\infty$ for every $j\in J$, therefore $i_{1}\in I$ and $B_{ii_1}(1)=A_{ii_1}(1)\neq -\infty$.

Let us assume that the sequence is defined up to rank $n$.
Since $A(n+1)$ has no line of~$-\infty$, there exists an $i_{n+1}\in [1,\cdots,d]$, such that $A_{i_ni_{n+1}}(n+1)\neq -\infty$. 

Since $\omega\in\A_i$, we have:
$$-\infty=\left(B(1)\cdots B(n)D(n+1)0\right)_i\ge \left(B(1)\cdots B(n)\right)_{ii_n}+\left(D(n+1)0\right)_{i_{n}},$$
therefore $\left(D(n+1)0\right)_{i_{n}}=-\infty$.

It means that every entry on  line $i_n$ of $D(n+1)$ is $-\infty$, that is \mbox{$A_{i_nj}(n+1)=-\infty$} for every $j\in J$, therefore $i_{n+1}\in I$ and \mbox{$B_{i_ni_{n+1}}(n+1)=A_{i_ni_{n+1}}(n+1)\neq -\infty$.}

Finally, we have:
$$\left(B(1)\cdots B(n+1)\right)_{ii_{n+1}}\ge\left(B(1)\cdots B(n)\right)_{ii_{n}}+B_{i_ni_{n+1}}(n+1)\neq -\infty .$$

\item To every matrix $A\in\rma^{d\times d}$, we associate the matrix $\widehat{A}$ defined by $\widehat{A}_{ij}=-\infty$ if $A_{ij}=-\infty$ and $A_{ij}=0$ otherwise. For every matrix $A,B\in\rma^{d\times d}$, we have $\widehat{A  B}=\widehat{A} \widehat{B}$.\\

The sequence defined by $R(n):=\widehat{B(1) \cdots  B(n)}$ is a Markov chain whose space of states is $\{0,-\infty\}^{I\times I}$ and whose transitions are defined by:
$$\P\left(R(n+1)=F|R(n)=E\right)=\P\left(\widehat{E  B(1)}=F\right).$$
For every $i,j\in I$, we have $R_{ij}(n)=0$ if and only if $\left(B(1) \cdots  B(n)\right)_{ij}\neq-\infty$.\\

Let $E$ be a recurrent state of this chain. Let us assume there exists a $k\in [1,\cdots,d]$ such that $E_{ik}=0$. Then, since $\G(B)$ is strongly connected, for every $j\in I$, there exists a $p\in\N$, such that $\left(B(1)\cdots B(p)\right)_{kj}\neq-\infty$ with positive probability, therefore there exists a state $F$ of the chain, reachable from state $E$ and such that $F_{ij}=0$. Since $E$ is recurrent, so is $F$ and $E$ and $F$ are in the same recurrence class.

Let us chose  $(i,l)\in I^2$ for a while. In each recurrence class of the Markov chain, either there exists a matrix $F$ such that $F_{il}=0$, or every matrix has only $-\infty$ on line~$i$.\\

Now, let us chose $(l,j)\in I\times J$, such that $\P(A_{lj}(1)\neq-\infty)>0$. Let $\mathcal{E}$ be a set with exactly one matrix $F$ in each recurrence class, such that $F_{il}=0$ whenever there is such a matrix in the class.
Let $S_n$ be the $n^\mathrm{th}$ time $\left(R(m)\right)_{m\in\N}$ is in $\mathcal{E}$.

Since the Markov chain has finitely many states and $\mathcal{E}$ intersect every recurrence class, $S_n$ is almost-surely finite. By definition of $\mathcal{E}$, we have almost-surely either $\left(B(1)\cdots B(S_n)\right)_{il}\neq-\infty$ or \mbox{$\forall m\in I, \left(B(1)\cdots B(S_n)\right)_{im}=-\infty$.}

It follows from~$i)$ that, if $\omega\in\A_i$, we are in the first situation. Therefore, we have, for every~$N\in\N$:
\begin{equation}\label{majprob}
\P\left[\A_i\right]\le\P\left[\forall n\in [1,\cdots,N], \left(D(S_n+1)0\right)_l=-\infty \right].
\end{equation} 

Conditioning the event $\left\{ \forall n\in [1,\cdots,N], \left(D(S_n+1)0\right)_l=-\infty \right\}$ by $S_N$, we have
\begin{eqnarray*}
\lefteqn{\P\left[\forall n\in [1,\cdots,N], \left(D(S_n+1)0\right)_l=-\infty \right]}\\
&=&\sum_{k\in\N}\P\left[S_N=k, \forall n\in [1,\cdots,N], \left(D(S_n+1)0\right)_l=-\infty \right]\\
&=&\sum_{k\in\N}\P\left[S_N=k, \left(D(k+1)0\right)_l=-\infty, \forall n\in [1,\cdots,N-1], \left(D(S_n+1)0\right)_l=-\infty \right]\\
&=&\sum_{k\in\N}\P\left[\left(D(k+1)0\right)_l=-\infty\right]\P\left[S_N=k, \forall n\in [1,\cdots,N-1], \left(D(S_n+1)0\right)_l=-\infty \right]\\
&=&\P\left[\left(D(1)0\right)_l=-\infty\right]\P\left[\forall n\in [1,\cdots,N-1], \left(D(S_n+1)0\right)_l=-\infty \right],
\end{eqnarray*}
because $\left\{\omega\in\Omega\left|S_N=k, \forall n\in [1,\cdots,N-1], \left(D(S_n+1)0\right)_l=-\infty \right.\right\}$ only depends on random matrices  $A(1),\cdots,A(k)$.

Finally, we have, for every $N\in\N$: 
$$\P\left[\forall n\in [1,\cdots,N], \left(D(S_n+1)0\right)_l=-\infty \right]=\left(\P\left[\left(D(1)0\right)_l=-\infty\right] \right)^N.$$
Because of the choice of $l$, $\P\left[\left(D(1)0\right)_l=-\infty\right] \le \P\left[A_{lj}(1)0=-\infty\right]<1$
and because of Equation~(\ref{majprob}), $\P\left(\A_i\right)=0.$
\end{enumerate} 
\end{proof}

\subsection{Fixed structure case}
Under the assumptions of Theorem~\ref{thSF1}, the hypotheses of Theorem~\ref{thgene} readily follows from the next two lemmas. Therefore Theorem~\ref{thSF1} is a consequence of Theorem~\ref{thgene}.
\begin{lem}\label{lemSF1}
Let $(\Omega,\theta,\P)$ be a measurable dynamical system such that for every $k\in [1,\cdots,d]$, $\theta^k$ is ergodic, and let $A:\Omega\rightarrow\rma^{d\times d}$ be a random matrix with no line of~$-\infty$, such that $\max_{ij} A_{ij}^+(0)$ is integrable.
If $A$ has fixed structure and $\G(A)$ is strongly connected, then the $y(n,0)$ associated to $A(n)=A\circ\theta^n$ satisfy
$$\lim_n\frac{1}{n}y(n,0)=\gamma(A)\1~\mathrm{ a.s.\,.}$$
\end{lem}

\begin{lem}\label{lemSF}
If a random matrix $A$ has fixed structure and has no line of~$-\infty$, then for every component $c_l$ of $\G(A)$, the random matrix $A^{\{l\}}$ has no line of~$-\infty$.
\end{lem}

\begin{proof}[Proof of Lemma~\ref{lemSF}]
Let $c_l$ be a component of $\G(A)$ and let us chose every $i\in H_l$.

If $i$ is in a component  $c_k$ of $\G(A)$, then there exists a path from $i$ to a component $c_m$ of $\G(A)$ such that $\gamma^{(m)}=\gamma^{[k]}$.
Let $j$ be the first node after $i$ on this path. Let $c_p$ be the component of $j$. Since $k\rightarrow p\rightarrow m$, we have: $$\gamma^{[k]}\ge\gamma^{[p]}\ge\gamma^{[m]}\ge\gamma^{(m)}=\gamma^{[k]}$$
and finally $\gamma^{[p]}=\gamma^{[k]}$. Because of proposition~\ref{propGiA}~$v)$, $\gamma^{[k]}= \gamma^{[l]}$, therefore  $\gamma^{[p]}= \gamma^{[l]}$, that is $p\in G_l$, and $j\in H_l$.

By definition of $\G(A)$, we have $\P(A_{ij}\neq-\infty)>0$, but because of fixed structure, it means $\P(A_{ij}\neq-\infty)=1$. Therefore $A^{\{l\}}$ has no line of~$-\infty$.
\end{proof}

We end this section with the proof of Lemma~\ref{lemSF1}\,.
\begin{proof}[Proof of Lemma~\ref{lemSF1}]
Let $L$ be the limit of $\frac{1}{n}y(n,0)$, which exists according to~\cite{BouschMairesse}.

Because of $\G(A)$'s strong connectivity and the fixed structure, for every entries $i,j$,  there exists $k_{ij}\in [1,\cdots,d]$, such that:
$$\left(A(-1)\cdots A(-k_{ij})\right)_{ij}\neq -\infty\textrm{ a.s.\,.}$$

It implies that
\begin{equation}\label{inegLiLj}
L_i\ge L_j\circ \theta^{-k_{ij}}~\mathrm{ a.s.\,.}
\end{equation} 

Practically, for $i=j$, it implies $L_i\ge L_i\circ\theta^{k_{ii}}$ almost-surely. Therefore $L_i= L_i\circ\theta^{k_{ii}}$ a.s.\,, and because of the ergodicity of $\theta^{k_{ii}}$, $L_i$ is almost-surely constant.

Equation~(\ref{inegLiLj}) therefore becomes $L_i\ge L_j$, and by symmetry $L_i=L_j$.  Finally we have for every $i\in [1,d]$:
$$L_i=\max_jL_j=\lim_n\frac{1}{n}\max_jy_j(n,0)=\gamma(A)\textrm{ a.s.\,.}$$
\end{proof}

\subsection{Precedence case}
In this section, we show that the hypotheses of Theorem~\ref{thprec} imply those of Theorem~\ref{thgene}. 
Hypothesis~$1.$ is obvious because of the precedence condition.
Hypotheses $2.$ and $3.$ both follow from the next lemma, whose proof is postponed to the end of the section:

\begin{lem}\label{lemprec}
Let $\sAn$ be a sequence satisfying the hypotheses of Theorem~\ref{thprec}. If $\G(A)$ is strongly connected, then for every $i\in [1,\cdots,d]$, there exists a random variable  $N$ with values in $\N$ such that for every $n\ge N$
\begin{equation}\label{eqlemprec}
\forall j\in [1,\cdots,d], \left(A(-1)\cdots A(-n)\right)_{ij}\neq -\infty
\end{equation}
\end{lem}

Let us assume the hypotheses of Theorem~\ref{thprec}.
Without loss of  generality, we assume that $A(n)=A\circ\theta^n$, where $(\Omega,\theta,\P)$ is a measurable dynamical system and $A$ is a random matrix.

To check hypothesis $2.$, we deduce from the lemma that if $\G(A)$ is strongly connected, then for every $i\in [1,\cdots,d]$, 
$$y_i(n,0)\ge \min_j\left(A(-1)\cdots A(-N)\right)_{ij}+\max_jy_j(n-N,0)\circ\theta^{-N},$$
and therefore
$$\liminf_n\frac{1}{n}y_i(n,0)\ge \lim_n \frac{1}{n} \min_j\left(A(-1)\cdots A(-N)\right)_{ij}+ \lim_n\frac{1}{n}\max_jy_j(n-N,0)\circ\theta^{-N} =\gamma(A).$$
Because of the definition of~$\gamma(A)$, we also have $\limsup_n\frac{1}{n}y_i(n,0)\le \gamma(A)$, therefore $\lim_n\frac{1}{n}y_i(n,0)= \gamma(A)$.

We apply this result to $A^{(m)}$ where $c_m$ is a dominating component, and this proves that Hypothesis~$2.$ is satisfied.\\

To check Hypothesis~$3.$, we apply  Lemma~\ref{lemprec} to matrix $B$ of decomposition~(\ref{decompblocs}), and we conclude the proof thanks to the ergodicity, that ensures there exists $n\ge N$ such that $D_{i'j}(-n-1)\neq-\infty$, provided $\P\left(D_{i'j}(1)\neq-\infty\right)$. Since there is such a pair $(i',j)$, it proves:
$$\left(B(-1)\cdots B(-n)D(-n-1)0\right)_i\ge \left(B(-1)\cdots B(-n)\right)_{ii'}+D_{i'j}(-n-1)>-\infty,$$
that is Hypothesis~$3.$ is checked.

\begin{proof}[Proof of Lemma~\ref{lemprec}]
Because of the ergodicity of $\theta$, if \mbox{$\P\left(A_{ij}(0)\neq-\infty\right)>0$,} then there exists almost-surely an $n_{ij}$ such that $A_{ij}(n_{ij})\neq-\infty$. That being the case, Poincar\'e recurrence theorem states that there are infinitely many such $n_{ij}$.\\

Let us chose $i\in [1,d]$.
Because of the precedence condition, the sequence of sets $$\A(n)=\left\{j\left|\left(A(-1)\cdots A(-n)\right)_{ij}\neq -\infty\right.\right\}$$ is increasing with $n$. Let us show that, from some rank $\A(n)=[1,d]$.\\

Since $\G(A)$ is strongly connected, there exists for every $j$ a finite sequence $i_0=i$, $i_1,\cdots, i_k=j$ such that for every $l\in [1,k]$, $\P\left(A_{i_li_{l+1}}(0)\neq-\infty\right)>0$. Because of what we said in the first paragraph of the proof, there are almost-surely $n_1<\cdots<n_k$ such that $A_{i_li_{l+1}}(n_{l})\neq-\infty$. Then, for every $n\ge n_i$, $\left(A(-1)\cdots A(-n)\right)_{ii_{l+1}}\neq-\infty$, and for every $n\ge n_k$, we have $\left(A(-1)\cdots A(-n)\right)_{ij}\neq-\infty$, that is $j\in \A(n)$. Since it is true for every $j$, there exists $N\in\N$ such that for every $n\ge N$, $\A(n)=[1,\cdots,d]$, which  concludes the proof of the lemma.

\section{Acknowledgements}
I have done this work both during my PhD at Universit\'e de Rennes~1 and as a JSPS postdoctoral fellow at Keio University. During this time, many exchange with J.~Mairesse have been a great help. This paper owes much to him.
\end{proof}
\bibliographystyle{alpha}
\bibliography{max+}

\def\cprime{$'$} \def\cprime{$'$}
\begin{thebibliography}{HOvdW06}

\bibitem[Bac92]{Baccelli}
F.~Baccelli.
\newblock Ergodic theory of stochastic petri networks.
\newblock {\em Annals of Probability}, 20(1):375--396, 1992.

\bibitem[BCOQ92]{BCOQ}
F.~Baccelli, G.~Cohen, G.J. Olsder, and J.P. Quadrat.
\newblock {\em Synchronisation and Linearity}.
\newblock John Wiley and Sons, 1992.

\bibitem[BH00]{TCPmp}
F.~Baccelli and D.~Hong.
\newblock Tcp is max-plus linear and what it tells us on its throughput.
\newblock In {\em SIGCOMM 00:Proceedings of the conference on Applications,
  Technologies, Architectures and Protocols for Computer Communication}, pages
  219--230. ACM Press, 2000.

\bibitem[BL92]{BaccelliLiu}
F.~Baccelli and Z.~Liu.
\newblock On a class of stochastic recursive sequences arising in queuing
  theory.
\newblock {\em The Annals of Probability}, 20(1):350--374, 1992.

\bibitem[BM]{BouschMairesse2}
T.~Bousch and J.~Mairesse.
\newblock Communication personnelle.

\bibitem[BM03]{BouschMairesse}
T.~Bousch and J.~Mairesse.
\newblock Fonctions topicales à portée finie et fonctions uniformément
  topicales.
\newblock Technical Report 2003-002, LIAFA (CNRS, UMR 7089 Université Paris 7),
  2003.
\newblock à paraitre dans \textit{Dynamical Systems}.

\bibitem[Bra93]{Braker}
H.~Braker.
\newblock {\em Algorithms and Applications in Timed Discrete Event Systems}.
\newblock PhD thesis, Delft University of Technology, Dec 1993.

\bibitem[CDQV85]{cohen85a}
G.~Cohen, D.~Dubois, J.P. Quadrat, and M.~Viot.
\newblock A linear system theoretic view of discrete event processes and its
  use for performance evaluation in manufacturing.
\newblock {\em IEEE Trans. on Automatic Control}, AC--30:210--220, 1985.

\bibitem[Coh88]{Cohen}
J.E. Cohen.
\newblock Subadditiviy, generalised products and operations research.
\newblock {\em SIAM Review}, 30(1):69--86, 1988.

\bibitem[dKHA03]{RailwayMpKHA}
A.~F. de~Kort, B.~Heidergott, and H.~Ayhan.
\newblock A probabilistic {$(\max,+)$} approach for determining railway
  infrastructure capacity.
\newblock {\em European J. Oper. Res.}, 148(3):644--661, 2003.

\bibitem[GM99]{GaubertMairesseIEEE}
S.~Gaubert and J.~Mairesse.
\newblock Modeling and analysis of timed {P}etri nets using heaps of pieces.
\newblock {\em IEEE Trans. Automat. Control}, 44(4):683--697, 1999.

\bibitem[Gri90]{Griffiths}
R.~B. Griffiths.
\newblock Frenkel-{K}ontorova models of commensurate-incommensurate phase
  transitions.
\newblock In {\em Fundamental problems in statistical mechanics VII (Altenberg,
  1989)}, pages 69--110. North-Holland, Amsterdam, 1990.

\bibitem[Hei00]{CaractMpQueuNet}
B.~Heidergott.
\newblock A characterisation of {$(\max,+)$}-linear queueing systems.
\newblock {\em Queueing Systems Theory Appl.}, 35(1-4):237--262, 2000.

\bibitem[Hon01]{Hong}
D.~Hong.
\newblock Lyapunov exponents: When the top joins the bottom.
\newblock Technical Report RR-4198, INRIA,
  http://www.inria.fr/rrrt/rr-4198.html, 2001.

\bibitem[HOvdW06]{MpAtWork}
B.~Heidergott, G.~J. Oldser, and J.~van~der Woude.
\newblock {\em Max plus at work}.
\newblock Princeton Series in Applied Mathematics. Princeton University Press,
  Princeton, NJ, 2006.
\newblock Modeling and analysis of synchronized systems: a course on max-plus
  algebra and its applications.

\bibitem[Mai95]{theseMairesse}
J.~Mairesse.
\newblock {\em Stabilité des systèmes à événements discrets stochastiques.
  Approche algébrique}.
\newblock PhD thesis, \'Ecole polytechnique, 1995.

\bibitem[Mai97]{Mairesse}
J.~Mairesse.
\newblock Products of irreducible random matrices in the $(\max,+)$ algebra.
\newblock {\em Adv. in Appl. Probab.}, 29(2):444--477, 1997.

\bibitem[Mer05]{theseGM}
G.~Merlet.
\newblock {\em Produits de matrices al\'eatoires~: exposants de Lyapunov pour
  des matrices al\'eatoires suivant une mesure de Gibbs, th\'eor\`emes limites
  pour des produits au sens max-plus.}
\newblock PhD thesis, Universit\'e de Rennes, 2005.
\newblock http://tel.ccsd.cnrs.fr/tel-00010813.

\bibitem[Vin97]{vincent}
J.-M. Vincent.
\newblock Some ergodic results on stochastic iterative discrete events systems.
\newblock {\em Discrete Event Dynamic Systems}, 7(2):209--232, 1997.

\end{thebibliography}

\end{document}